\documentclass[11pt,twoside]{article}
\include{epsf}
\usepackage{latexsym, amssymb, amsmath, theorem}
\usepackage{calc, epsfig}
\numberwithin{equation}{section}

\theoremstyle{change}
\theorembodyfont{\itshape}
\newtheorem{theorem}{Theorem}[section]
\newtheorem{proposition}[theorem]{Proposition}
\newtheorem{lemma}[theorem]{Lemma}
\newtheorem{corollary}[theorem]{Corollary}

\theorembodyfont{\rmfamily}
\newtheorem{definition}[theorem]{Definition}

\newtheorem{remark}[theorem]{Remark}
\newenvironment{proof}{{\noindent \textbf{Proof}\,\,}}{\hspace*{\fill}$\Box$\medskip}
\def\rd{\mathbb R^2}
\def\cc{\mathbb C}
\def\bc{\overline{\cc}}
\def\ci{C^{\infty}}

\title{\textbf 
{Simple proofs of uniformization theorems}}
\author{A.A.Glutsyuk\\[5pt]
\textit{CNRS, Unit\'e de Math\'ematiques Pures et Appliqu\'ees, M.R.,} \\
\textit{\'Ecole Normale Sup\'erieure de Lyon} \\
\textit{46 all\'ee d'Italie, 69364 Lyon Cedex 07 France}}
\begin{document}
\maketitle
\def\td{\mathbb T^2}
\begin{abstract} 
The measurable Riemann mapping theorem proved 
by Morrey and in some particular cases by Ahlfors, Lavrentiev 
and Vekua, says that any measurable almost complex 
structure on $\rd$ ($S^2$) with bounded dilatation is integrable: there is a 
quasiconformal homeomorphism of $\rd$ ($S^2$) onto $\cc$ 
($\bc$) transforming the given almost complex structure to the standard one. 
We give an elementary proof of this theorem that is done as follows. 
Firstly we prove its double-periodic version: each $\ci$ almost complex 
structures on the two-torus can be transformed by a diffeomorphism 
 to the standard complex structure on appropriate complex torus. The proof is 
 based on the homotopy method for the Beltrami equation on $\td$ 
with parameter. (As a by-product, 
we present a simple proof of the Poincar\'e-K\"obe theorem saying that each 
simply-connected Riemann surface is conformally equivalent to either 
$\overline{\cc}$, or $\cc$, or the unit disc.) Afterwards the 
general case is treated by $\ci$ double-periodic approximation  and simple
 normality arguments (involving Gr\"otzsch inequality) following 
 the classical scheme.  
\end{abstract}

\tableofcontents
\vspace{.5cm}
\hrule

\section{Introduction, the plan of the paper and history}
\subsection{Uniformization theorems. The plan of the paper}
A linear complex structure on $\rd$ is a structure of a linear space 
over $\mathbb C$ (we fix an orientation and consider it to be 
compatible with the complex structure). 
The {\it (almost) complex structure} on a real two-dimensional surface is 
a family of linear complex structures on the tangent planes at its points.  
A linear complex structure on $\rd$ defines an ellipse 
in $\rd$ centered at 0, which is an orbit under the $S^1$- action 
by multiplication by complex numbers with unit module. (The ellipse 
corresponding to the standard complex structure on $\cc$ is a circle.) 
The {\it dilatation} of a nonstandard linear complex structure on $\cc$ 
(with respect to the standard complex structure) is the excentricity 
of the corresponding ellipse (i.e., the ratio of the largest radius over 
the smallest one). An almost complex 
structure defines an ellipse field in tangent planes, and vice versa: 
the ellipse field determines the almost complex structure in a unique way. 

If our surface is a Riemann surface (with a fixed complex structure), 
then any  (nonstandard) almost complex 
structure has a well-defined dilatation at each point of the surface. 
In this case an almost complex structure is said to be {\it bounded}, 
if its dilatation is bounded. 
The  {\it (total) dilatation} of a bounded almost complex structure is 
the supremum of its dilatations (more precisely, the minimal supremum 
of dilatations after possible correction of the almost complex 
structure over a measure zero set). 

Each real linear isomorphism $\cc\to\cc$ acts on the space of the ellipses centered 
at 0, and hence, on the space of linear complex structures. Its {\it 
dilatation} is defined to be the dilatation of the image of the 
standard complex structure (which is equal to the excentricity of the 
image of a circle centered at 0). 
The action of a differentiable homeomorphism of domains in $\mathbb C$ 
on the almost complex structures 
and its dilatation (at a point) are defined to be those of its derivative. 
Its {\it (total)} dilatation is the supremum of the dilatations through all 
the points. 

It appears that {\it any $\ci$ (and even measurable) bounded almost 
complex structure is integrable}, that is, can be transformed to a true  
complex structure by a $\ci$ (respectively, quasiconformal) homeomorphism, see 
the following Definition and Theorem. 

\begin{definition} (see, e.g., [Ah2]). Let $K>0$. 
A homeomorphism of domains in $\cc$ is said to be $K$- 
{\it quasiconformal} (or $K$- homeomorphism), 
if it  has local $L_2$ derivatives and its dilatation 
(at the differentiability points with nonzero derivative) is 
 no greater than $K$. 
A homeomorphism is said to be quasiconformal if it is 
$K$- quasiconformal for some $K>0$. 
\end{definition}

\begin{remark} \label{grr} The dilatations of a 
differentiable homeomorphism and its inverse are equal.  In particular, the inverse to a $K$- diffeomorphism is also a 
$K$- diffeomorphism. The composition of two $K$- diffeomorphisms is a 
$K^2$- diffeomorphism. This follows from definition. 
\end{remark}

\begin{proposition} \label{group}  (see [Ah2])   
The quasiconformal homeomorphisms of a Riemann surface form a group.
\end{proposition}

\begin{proposition} \label{measure0} The image of a zero measure 
set under a quasiconformal homeomorphism has also zero measure.  
\end{proposition}

\begin{corollary} \label{cormeasure} For any quasiconformal homeomorphism 
the set of its differentiability points with zero derivative has measure 
zero. 
\end{corollary} 

\begin{proof} The image of the set from the Corollary has zero measure 
by definition. Therefore, the set itself has zero measure (Proposition 
\ref{measure0} applied to the inverse mapping, which is quasiconformal 
by Proposition \ref{group}).
\end{proof}  
 
Both Propositions are proved in Subsection 3.4 and neither them, nor the 
Corollary  will be used in the paper. 

\begin{definition}
A homeomorphism $\cc\to\cc$ is said to be {\it normalized}, if it fixes 0 and 1. 
\end{definition}

\begin{theorem} \label{c} ([AhB], [M]). For any measurable 
bounded almost complex structure $\sigma$ on $\mathbb C$ there exists a unique 
normalized quasiconformal homeomorphism 
$\mathbb C\to\mathbb C$ that transforms $\sigma$ 
to the standard complex structure 
(at the differentiability points with nonzero derivative). 
If $\sigma$ is $\ci$ in some domain, then 
the homeomorphism is a $\ci$ diffeomorphism while restricted to this domain.  
\end{theorem}

{\bf Addendum [AhB].} If a bounded almost
complex structure on $\cc$ varies analytically in a complex parameter, then 
so does the corresponding homeomorphism from Theorem \ref{c}.  
\medskip

\begin{remark} A quasiconformal homeomorphism of a once punctured domain 
extends quasiconformally to the puncture 
(in particular, the homeomorphism from Theorem \ref{c} is quasiconformal 
at infinity). This follows easily 
from the local uniqueness of the quasiconformal homeomorphism up to 
composition with conformal mapping (Proposition \ref{uniq}, see 
Subsection 3.4) and the theorem on erasing isolated singularities 
of bounded holomorphic functions. 
\end{remark}

In the present paper we give proofs of Theorem \ref{c} (Sections 2, 3) 
and the Addendum (Subsection 3.5) 
that seem to be simpler than the known proofs and easier to explain. 
A historical overview will be given in Subsection 1.4. 

\begin{remark} The proof of the local integrability of an analytic almost 
complex structure is elementary: it is done immediately 
by analyzing the complexification of the corresponding $\cc$- linear 1- form 
(this proof is due to Gauss). But it is already nontrivial in the $\ci$ case. 
\end{remark}

The measurable versions of the Theorem and the Addendum have many 
very important 
applications in various domains of mathematics, especially in holomorphic 
dynamics and the 
Kleinian group theory (quasiconformal surgery, where one deals with 
invariant almost complex structures that are discontinuous...), see, e.g.,   
[CG]. 

For the proof of Theorem \ref{c} we firstly prove (in Section 2) its 
version for $\ci$ almost 
complex structures on the two-torus: the proof uses only elementary 
Fourier analysis. 

\begin{theorem} \label{td} ([Ab]) For any $\ci$ 
bounded almost complex structure $\sigma$ on $\td$ there exists a 
$\ci$ diffeomorphism of $\td$ onto appropriate 
complex torus (the latter torus depends on $\sigma$) 
that transforms $\sigma$ to the standard complex structure.
\end{theorem}

Then in Section 3 we deduce Theorem \ref{c} from 
Theorem \ref{td} by using double-periodic approximations of a given  
almost complex structure on $\cc$ and simple normality arguments  involving a 
Gr\"otzsch inequality for annuli diffeomorphisms. This deduction 
follows the classical scheme [Ah2]. 

The proof of  Theorem \ref{td} presented below is implicitly contained in the 
previous paper [Gl] by the author, where the same method was used to 
prove a foliated version of Theorem \ref{c}. We prove the existence of a global 
nowhere vanishing $\sigma$- holomorphic differential. To do this, 
we use the homotopy method for the Beltrami equation with parameter, which   
reduces the proof to solving a linear 
ordinary differential equation in $L_2(\td)$. We prove
regularity of its solution by 
showing that the equation is bounded in any Sobolev space $H^s(\td)$.

In Subsection 1.3 we give a proof of  
the classical Poincar\'e-K\"obe uniformization theorem using 
Theorem \ref{td}: 
\begin{theorem} \label{pk} [Ko1], [Ko2], [P]. Each simply-connected Riemann surface is 
conformally equivalent to either unit disc, or $\cc$, or the Riemann sphere.
\end{theorem}

In the proofs of the previously mentioned Theorems we use  the 
well-known notations (recalled in the next Subsection) concerning 
almost complex structures. 
\def\omu{\omega_{\mu}}

\subsection{Complex structures and uniformizing differentials. Basic notations}

 To a (nonstandard) almost complex structure (denoted $\sigma$) 
 on a subset $D\subset \mathbb C$ 
we put into correspondence a $\cc$- valued 1- form that is 
$\cc$- linear with respect 
to $\sigma$. The latter form can be normalized to have the type 
\begin{equation}
\omu=dz+\mu(z) d\bar z, \ |\mu|<1.\label{1.1}
\end{equation}
 The function $\mu$ is uniquely defined by $\sigma$. Vice versa, for arbitrary 
complex-valued function $\mu$, $|\mu|<1$, the 1- form (\ref{1.1}) defines the  
unique complex structure for which it is $\cc$- linear. We denote by 
$\sigma_{\mu}$ the almost complex structure thus defined (whenever the contrary 
is not specified). Then $\sigma_{\mu}$ is bounded, if and only if $\sup|\mu|<1$. 

\begin{remark} The ellipse associated to $\sigma_{\mu}$ 
on the tangent plane at a point $z$ is given 
by the equation $|dz+\mu(z)d\bar z|=1$; the dilatation (excentricity) is 
equal to $\frac{1+|\mu(z)|}{1-|\mu(z)|}$. 
\end{remark}
We will be looking for a 
differentiable homeomorphism $\Phi(z)$ that is holomorphic, i.e., 
that transforms $\sigma_{\mu}$ to the standard complex structure. This is equivalent 
to say that the differential of $\Phi$ (which is a closed form) 
is a $\cc$- linear form, i.e., has the type $f(z)(dz+\mu d\bar z)$:  
$$\frac{\partial\Phi}{\partial\bar z}=\mu\frac{\partial\Phi}{\partial z}.$$

\begin{remark} Conversely, let $\mu$ be $\ci$, $|\mu|<1$. 
Then any $\ci$ closed 
1- form $f(z)(dz+\mu d\bar z)$ is $\sigma_{\mu}$-  holomorphic, 
i.e., is a differential of a complex-valued $\ci$ function $\Phi$ transforming 
$\sigma_{\mu}$ to the standard complex structure. 
A form $f(z)(dz+\mu d\bar z)$ is closed if and only if 
\begin{equation}\partial_{\bar z}f=\partial_z(\mu f).\label{dz}\end{equation}
\end{remark}
\begin{definition} A Riemann surface is said to be {\it parabolic}, 
if its universal covering is conformally equivalent to $\mathbb C$ 
(i.e, the surface is either $\mathbb C$, or $\mathbb C^*$, or a complex torus). 
\end{definition}

\begin{definition} \label{unifdif} The {\it uniformizing differential} on 
$\mathbb C$ (or on a complex torus) 
with the affine coordinate $z$ is the 1- form $dz$ or its nonzero constant 
multiple. More generally, a holomorphic 1- form on a parabolic Riemann surface 
is said to be {\it a uniformizing differential}, if the primitive of 
its lifting to the universal cover is a conformal isomorphism onto $\mathbb C$. 
\end{definition}

\def\smu{\sigma_{\mu}}

\begin{remark} \label{umetric} The uniformizing differential is well-defined up 
to multiplication by constant. It coincides with the unique (up to constant) 
nowhere vanishing holomorphic differential whose squared module is 
a complete metric. 
\end{remark}
\begin{proposition} \label{uniftd} Let $\mu:\td\to\cc$ be a $\ci$ function, 
$|\mu|<1$.  Suppose there is a $\ci$ 
nowhere vanishing function $f:\td\to\cc\setminus0$ satisfying (\ref{dz}). 
Then the corresponding almost complex structure $\sigma_{\mu}$ is 
integrable and the form $f\omu$ 
 is a uniformizing differential of $(\td,\smu)$. 
\end{proposition}

The Proposition follows from compactness and the two previous Remarks.

\subsection{Proof of the uniformization Theorem \ref{pk} modulo 
Theorem \ref{td}}

Let $S$ be a simply-connected Riemann surface. Then it is either contractible, 
then homeomorphic $\rd$, or is homeomorphic to the two-sphere. 
We prove the statement of Theorem \ref{pk} only in 
the case, when $S$ is contractible: we show that 
$S$ is conformally equivalent to 
either $\cc$ or disc. Then if $S$ is sphere, it follows that 
$S$ is conformally-equivalent to $\overline{\cc}$ (by the previous 
statement applied to once punctured $S$ and the theorem on 
erasing isolated singularities of bounded holomorphic functions). 
In the proof of Theorem \ref{pk} we use the following Corollary 
of Theorem \ref{td} and 
the Riemann mapping theorem (saying that any simply-connected domain in $\cc$ 
distinct from $\cc$ is conformally equivalent to unit disc: the proof 
 is elementary and is contained in standard courses of complex analysis.)

\begin{corollary} For any bounded $\ci$ almost complex structure $\sigma$ 
on the closed unit disc $\overline D$ there exists a $\ci$ diffeomorphism of the open 
disc $D$ onto itself transforming $\sigma$ to the standard complex structure. 
\end{corollary}

\begin{proof} Let us extend $\sigma$ 
to $\rd$ up to a double-periodic bounded $\ci$ almost complex structure 
(say, with periods $4$ and $4i$) and consider the quotient torus equipped 
with the induced almost complex structure.  
Then the corresponding tori diffeomorphism from 
Theorem \ref{td} transforms the latter structure to the standard one. Its  
lifting to the universal covers 
transforms $D$ to a simply-connected domain in $\cc$ and 
sends $\sigma$ to the standard complex structure. Now applying 
the Riemann mapping theorem to the image of $D$ proves the Corollary. 
\end{proof}

We assume that the Riemann surface $S$ is contractible, hence, 
admits a $\ci$ 1-to-1 parametrization by $\rd$. 
Its complex structure induces a $\ci$ almost complex structure 
(denote it $\sigma$) on $\rd$  
(not necessarily bounded). Take a growing sequence of discs  
$S_1\Subset S_2\Subset\dots S$ exhausting $S$ centered at 0. 
On each $S_n$ the almost complex structure $\sigma$ is bounded. By the 
Corollary, for any $n$ 
there is a diffeomorphism $\phi_n:S_n\to D$ conformal with respect 
to the complex structure of $S$, $\phi_n(0)=0$. 
Let $w$ be a local holomorphic chart on $S$ near 0, $w(0)=0$. 
Let us change $\phi_n$ to its constant multiple 
$\Phi_n=\lambda_n\phi_n$ having unit derivative in $w$ at 0. The family 
$\Phi_n$ is normal: each subsequence contains a subsequence converging 
uniformly on compact sets in $S$. Indeed, fix a $k\in\mathbb N$ 
and consider the $\ci$ injections 
$\Phi_n\circ\phi_k^{-1}:D\to\Phi_n(S_n)$, $n\geq k$. By construction, the 
latters are holomorphic and univalent, they 
send 0 to 0 and have one and the same derivative at 0. Therefore, they form a 
normal family, see [CG], hence, so do the $\Phi_n$'s. By construction, the 
limit of a converging subsequence of the $\Phi_n$' s is a conformal 
diffeomorphism of $S$ onto either a disc, or $\cc$. Theorem \ref{pk} is proved. 

\def\onu{\omega_{\nu}}

\subsection{Historical overview} The local integrability of a 
$\ci$ (and even H\"older) almost complex structure was proved by 
Korn [Korn] and Lichtenstein [Licht]; a simpler proof was obtained 
by Chern [Chern] and Bers [Be]. The local integrability together with the 
Poincar\'e-K\"obe uniformization Theorem \ref{pk} imply the  
global integrability  statement of Theorem \ref{c}. Lavrentiev [La] 
gave a direct proof of Theorem \ref{c} for continuous almost complex 
structures. Later Ahlfors [Ah1] and Vekua [Vek] gave another direct proofs  
under the previous (stronger) H\"older condition. 

In the general measurable case Theorem \ref{c} was proved by Morrey [M]. 
Later new proofs were obtained by Ahlfors and Bers [AhB], Bers and 
Nirenberg [BeN] and Boyarskii [Bo].  
(In fact, Lavrentiev and Morrey stated 
their theorems for almost complex structures on a disc, but their versions 
on $\rd$ follow immediately, e.g.,  by the arguments from the previous 
Subsection.) A new simpler proof of Theorem \ref{c}  
using $L_2$ analysis and Fourier transformation on $\rd$ 
was recently obtained by A.Douady and X.Buff [DB].  

\section{Smooth complex structures on $\td$. Proof of Theorem \ref{td}}

\subsection{Homotopy method. The sketch of the proof of Theorem \ref{td}}

Let $\mu:\td\to\cc$ be a $\ci$ complex-valued function, $|\mu|<1$, 
$\smu$ be the corresponding almost complex structure, see (\ref{1.1}). 
Theorem \ref{td} says that there exists 
a diffeomorphism transforming $(\td,\smu)$ into a complex torus equipped 
with the standard complex structure. To prove this statement, 
it suffices to construct 
a uniformizing differential, more precisely, a $\ci$ nowhere vanishing 
function $f:\td\to\cc\setminus0$ such that the form $f\omu$ is closed 
(see Proposition \ref{uniftd}), i.e., to solve partial differential 
equation (\ref{dz}) in a $\ci$ nowhere vanishing function $f$.

To solve (\ref{dz}), we use the homotopy method. 
Namely, we include $\smu$ into the one-parametric 
family of complex structures (denoted by $\sigma_{\nu}$) defined by  their 
$\mathbb C$- linear  1- forms 
$$\onu=dz+\nu(z,t)d\bar z,\ \nu(z,t)=t\mu(z),\ t\in[0,1].$$
The complex structure corresponding to the parameter value $t=0$ is the 
standard one, the given structure $\sigma_{\mu}$ corresponds to 
$t=1$. 
We will find a $C^{\infty}$ family $f(z,t):\td\times[0,1]\to\cc\setminus0$ of 
complex-valued nowhere vanishing $C^{\infty}$ functions on $\td$ depending on 
the same parameter $t$, $f(z,0)\equiv1$, such that the differential forms 
$f(z,t)\onu$ are closed, i.e., 
\begin{equation} \partial_{\bar z}f=\partial_z(f\nu).\label{dzt}\end{equation}
Then the function $f=f(z,1)$ is the one we are looking for. 

To construct the previous family of functions, we will find firstly 
a family $f(z,t)$ of 
{\it nonidentically-vanishing} (not necessarily nowhere vanishing) 
functions satisfying (\ref{dzt}): 

\begin{lemma} \label{homot} Let $\nu(z,t):\td\times[0,1]\to\cc$ be a $\ci$ family of 
$C^{\infty}$ functions on $\td$, $|\nu|<1$, $\nu(z,0)\equiv0$, $z$ be 
the complex coordinate on $\td$. There exists a $\ci$ family 
$f(z,t):\td\times[0,1]\to\cc$ of $\ci$ functions on $\td$ that are solutions of 
(\ref{dzt}) with the initial condition 
$f(z,0)\equiv1$ such that for any fixed $t\in[0,1]$ $f(z,t)\not\equiv0$ in $z$.
\end{lemma}

The Lemma will be proved in the next Subsection.
 
Below we show that in fact, the functions $f(z,t)$ from the Lemma vanish 
nowhere. To do this (and only in this place) we use the 
local integrability of  a $\ci$ complex structure:  

\begin{proposition} \label{loc} ([Korn], [Licht], [La], [Chern], [Be]). 
Let $D\subset\cc$ be a disc centered at 0, 
$\mu:D\to\cc$, $\mu\in\ci$, $|\mu|<1$, $\smu$ be the corresponding 
almost complex structure, see (\ref{1.1}). There exists a 
local $\sigma_{\mu}$- holomorphic univalent coordinate near 0.
\end{proposition} 

The Proposition will be proved in Subsection 2.3. 

\begin{proof} {\bf of Theorem \ref{td} modulo Lemma \ref{homot} 
and Proposition \ref{loc}.} Let $f(z,t)$ be a family of functions 
from the previous Lemma. By the previous discussion, it suffices 
to show that $f(z,t)\neq0$. This inequality holds for $t=0$, where $f=1$. 

 Let us prove that $f(z,t)\neq0$ by contradiction. Suppose the contrary. 
Then the set of the parameter values $t$ corresponding to the functions $f(z,t)$ 
having zeroes is nonempty (denote this set by $M$). Its complement 
$[0,1]\setminus M$ is open by definition. Let us show that the set $M$ is open 
as well. This will imply that the parameter segment is a union of two disjoint 
open sets, which will bring us to contradiction. It suffices to show that 
the (local) presense of a zero of a function $f$ persists under perturbation. 

Suppose $f(z_0,t)=0$ for some $z_0$, $t$ (let us fix them). It suffices to 
show that for $t'$ close to $t$ the function $f(z,t')$ has a zero 
near $z_0$. Let $w$ 
be the local holomorphic coordinate on $\td$ near $z_0$ from the previous 
Proposition corresponding to $\mu=\nu(z,t)$, $w(z_0)=0$. 
Suppose that 
the function $f(z,t)$ does not vanish identically on $\td$ 
 locally near $z_0$: one can achieve this by changing  
$z_0$, since $f$ does not vanish identically. Recall that $f\onu$ is 
a closed $\cc$- linear 1-form with  respect to the variable complex structure 
$\sigma_{\nu}$, hence, it is holomorphic in the 
coordinate $w$. Therefore, $f\onu=(w^k+\text{higher terms})dw$, $k\geq1$. 
Now by the index argument, the local presense of zero of $f$ on 
$\td$ persists under 
perturbation. This together with the previous discussion proves the inequality 
 $f(z,t)\neq0$ and Theorem \ref{td}.
 \end{proof}
 
 \subsection{Variable holomorphic differential: proof of Lemma \ref{homot}}
 
 Differentiating (\ref{dzt}) in $t$ yields (we denote 
 $\dot f$ the partial derivative in $t$ of a function $f$)
 \begin{equation}
\partial_{\bar z}\dot f-(\partial_z\circ\nu)\dot f=(\partial_z\circ\dot\nu)f.\label{2.2}
\end{equation}
where $\partial_z\circ\nu$ ($\partial_z\circ\dot\nu$) is the composition of the operator 
of the multiplication by the function $\nu$ (respectively, $\dot\nu$) and the 
operator $\partial_z$. 
Any solution $f$ of equation (\ref{2.2}) with the initial condition 
$f(z,0)\equiv1$ 
that vanishes identically on the torus for no value of $t$ is a one we 
are looking for. Let us show that (\ref{2.2}) is implied by a bounded 
linear differential equation in $L_2(\td)$. To do this, 
we use the following properties of the operators $\partial_z$ and $\partial_{\bar z}$. 
 
 \begin{remark} Denote $z=x_1+ix_2$, $x=(x_1,x_2)\in\rd$. 
 The operators $\partial_z$, $\partial_{\bar z}$ on $\td$ have common eigenfunctions 
 $e_n(x)=e^{i(n,x)}$, $n=(n_1,n_2)\in\mathbb Z^2$. The corresponding eigenvalues 
 (denote them $\lambda_n$ and $\lambda_n'$ respectively) have equal modules, 
 more precisely, 
 \begin{equation}\lambda_n'=-\overline{\lambda_n}.\label{2.3}\end{equation}
 This is implied by the fact that the operator $\partial_{\bar z}$ is conjugated to 
 $-\partial_z$ in the $L_2$ scalar product, which follows from definition. In fact, 
 $$\lambda_n=\frac i2(n_1-in_2),\ \lambda_n'=\frac i2(n_1+in_2).$$
 \end{remark}
 
 \begin{corollary} There exists a unique unitary operator $U:L_2(\td)\to 
 L_2(\td)$ preserving averages such that "$U=\partial_{\bar z}^{-1}\circ\partial_z$" 
 (more precisely, $U\circ\partial_{\bar z}=\partial_{\bar z}\circ U=\partial_z$). The operator 
 $U$ commutes with partial differentiations and extends up to a unitary operator 
 to any Hilbert Sobolev space of functions on $\td$. In particular, it preserves 
 the space of $\ci$ functions.
 \end{corollary}
 \begin{proof} 
 The operator $U$ from the Corollary is defined to have the previous 
 eigenfunctions $e_n$ with the eigenvalues $\frac{\lambda_n}{\lambda_n'}=
 \frac{n_1-in_2}{n_1+in_2}$. Its uniqueness follows immediately from the 
 previous operator equation on $U$ applied to the functions $e_n$. 
 The rest of the statements of the Corollary follow immediately 
 from definition and Sobolev embedding theorem (see [Ch], p.411). 
 \end{proof}
 
 Let us write down equation (\ref{2.2}) in terms of the new operator $U$. Applying 
the "operator" $\partial_{\bar z}^{-1}$ to (\ref{2.2}) and substituting 
$U=\partial_{\bar z}^{-1}\circ \partial_z$ yields  
$$(Id-U\circ\nu)\dot f=(U\circ\dot\nu)f.$$ 
This equation implies (\ref{2.2}). For any $t\in[0,1]$ the operator 
$Id-U\circ\nu$ in the left-hand side is invertible in $L_2(\td)$ and 
the norm of the 
inverse operator is bounded uniformly in $t$, since $U$ is unitary and the 
module  
$|\nu|$ is less than 1 and bounded away from 1 by compactness. 
Thus, the last equation can be rewritten as 
\begin{equation}
\dot f=(Id-U\circ\nu)^{-1}(U\circ\dot\nu)f,\label{2.6}
\end{equation}
 which is an ordinary differential equation in $f\in L_2(\td)$ with a 
 uniformly $L_2$- bounded operator in the right-hand side. 
  As it is shown below 
(in Proposition \ref{2.70}), the inverse $(Id-U\circ\nu)^{-1}$ is also uniformly 
bounded in each Hilbert Sobolev space $H^j(\td)$.  
 Therefore, equation (\ref{2.6}) written in arbitrary Hilbert Sobolev space has a 
 unique solution 
 with a given initial condition, in particular, with $f(z,0)\equiv1$ (the 
 theorem on existence and uniqueness of solution of ordinary 
 differential equation in Banach space with the right-hand side having 
 uniformly bounded derivative [Ch]). For any $t\in[0,1]$ this solution does 
 not vanish identically on $\td$ (uniqueness of solution) and belongs to 
 all the spaces $H^j(\td)$; hence, it is $C^{\infty}(\td)$ by Sobolev embedding 
 theorem (see [Ch], p.411). Thus, Lemma \ref{homot} is implied by the following
 
 \begin{proposition} \label{2.70} Let $x=(x_1,x_2)$ be affine coordinates on 
 $\mathbb R^2$, 
 $\td=\mathbb R^2\slash2\pi\mathbb Z^2$. 
 Let $s\geq0$, $s\in\mathbb Z$, $U$ be a linear operator 
 in the space of $C^{\infty}$ functions on $\td$ that commutes with the operators 
 $\frac{\partial}{\partial x_i}$, $i=1,2$, and extends to any Sobolev space 
 $H^j=H^j(\td)$, $0\leq j\leq s$, up to a unitary operator. Let $0<\delta<1$, 
 $\nu\in C^s(\td)$ be a complex-valued function, 
 $|\nu|\leq\delta$. The operator $Id-U\circ\nu$ is invertible and the 
 inverse operator 
 is bounded in all the spaces $H^j$, $0\leq j\leq s$. 
For any $0<\delta<1$, $j\leq s$, there 
 exists a constant $C>0$ (depending only on $\delta$ and $s$) such that for any 
 complex-valued  function $\nu\in C^s(\td)$ with $|\nu|\leq\delta$  
 $$||(Id-U\circ\nu)^{-1}||_{H^j}\leq C(1+\sum_{k\leq j}\max
 |\frac{\partial^k\nu}{\partial x_{i_1}, \dots,\partial x_{i_k}}|^j).$$
 \end{proposition}
 
 \begin{proof} Let us prove the Proposition for $s=1$. For higher $s$ its proof is 
 analogous. 
 \begin{equation}\text{By definition,}\ \ 
||U\circ\nu||_{L_2}\leq\delta<1.\label{2.7}
\end{equation}
 Hence, the operator 
 $Id-U\circ\nu$ is invertible in $L_2=H^0$ and 
\begin{equation}
(Id-U\circ\nu)^{-1}=Id+\sum_{k=1}^{\infty}(U\circ\nu)^k:\label{2.8}
\end{equation}
 the sum of the $L_2$ operator norms of the sum entries in 
  (\ref{2.8}) is  finite by (\ref{2.7}). 
 Let us show that the operator in the right-hand side of (\ref{2.8}) is well-defined 
 and bounded in $H^1$. 
 To do this, it suffices to show that the sum of the operator $H^1$- norms of 
 the same entries is finite. 
 
 Let $f\in H^1(\td)$. Let us estimate 
 $||(U\circ\nu)^kf||_{H^1}$. We show that for any $k\in\mathbb N$ 
\begin{equation}
||\frac{\partial}{\partial x_r}((U\circ\nu)^kf)||_{L_2}<ck\delta^{k-1}
||f||_{H^1}, \ \ c=\delta+\max|\frac{\partial\nu}{\partial x_r}|,\ \ 
r=1,2.\label{2.9}
\end{equation}

 This will imply the finiteness of the operator $H^1$- norm of the sum in the 
 right-hand 
 side of (\ref{2.8}) and Proposition \ref{2.70} 
 (with $C=4\sum_{k\in\mathbb N}k\delta^{k-1}=
 \frac4{(1-\delta)^2}$).  

Let us prove (\ref{2.9}), e.g., for $r=1$. The derivative in the 
left-hand side of (\ref{2.9}) equals  
$$(U\circ\nu)^k\frac{\partial f}{\partial
x_1}+\sum_{i=1}^k(U\circ\nu)^{k-i}
\circ(U\circ\frac{\partial\nu}{\partial x_1})\circ(U\circ\nu)^{i-1}f$$
(since $U$ commutes with the partial differentiation by the condition of 
Proposition \ref{2.70}). The $L_2$- norm of the first term in the previous formula 
is no greater than 
$\delta^k||f||_{H^1}$ by (\ref{2.7}).  Each term in its sum has $L_2$- norm 
no greater 
than $\delta^{k-1}\max|\frac{\partial\nu}{\partial x_1}|||f||_{L_2}$ by (\ref{2.7}). 
This proves (\ref{2.9}). The Proposition is proved. Lemma \ref{homot} 
is proved. 
\end{proof}
\begin{remark} The solution of equation (\ref{2.6}) with the initial condition 
 $f|_{t=0}\equiv1$ admits the following formula: 
 \begin{equation}
 f(x,t)=(Id-U\circ\nu)^{-1}(1)=
 1+U(\nu)+(U\circ\nu\circ U)(\dot\nu)+\dots\label{2.10}
 \end{equation}
 Indeed, its right-hand side is a well defined $\ci$  
 family of $C^{\infty}$ functions on $\td$, which follows from the 
 uniform boundedness of the operators $(Id-U\circ\nu)^{-1}$ in any given 
 Hilbert Sobolev space. By definition, it satisfies  the unit initial 
 condition. Differentiating (\ref{2.10}) in $t$ yields
 $$(Id-U\circ\nu)^{-1}\circ (U\circ\dot\nu)\circ(Id-U\circ\nu)^{-1}(1)=
(Id-U\circ\nu)^{-1}\circ(U\circ\dot\nu)f(x,t).$$
Hence, the function (\ref{2.10}) satisfies (\ref{2.6}).  
 \end{remark}
 
 \def\wt#1{\widetilde#1}
 \def\tdv{\mathbb T^2}

\subsection{Zero of holomorphic differential. Proof of Proposition \ref{loc}} 
 Let us prove the existence of local holomorphic coordinate. Without loss of 
 generality we assume that $\mu(0)=0$ (applying a linear change of variables).  
 One can achieve also that $\mu$ is arbitrarily 
 small with derivatives of orders up to 3 applying a homothety and taking 
 the restriction to a smaller disc centered at 0. 
 We consider that the disc where $\mu$ is defined 
is embedded into $\td$ and extend the 
 function $\mu$ smoothly to $\td$. We assume that the extended function 
 satisfies the inequality $||\mu||_{C^3(\td)}<\delta$; one can make $\delta$ 
 arbitrarily small. 
 
 Let $\nu(x,t)=t\mu$, $f(x,t)$ be the corresponding function 
 family from Lemma \ref{homot} constructed as the solution of 
differential equation (\ref{2.6}) with unit initial condition,  $f(x)=f(x,1)$. 
We show in the next paragraph 
 that $f(0)\neq0$, if the previous constant $\delta$ is small enough. 
 Then the local coordinate we are looking for is the function 
 $$w(z)=\int_0^zf(dz+\mu d\bar z).$$
 Indeed, it is well-defined and holomorphic by definition. Its local univalence 
 follows from the nondegeneracy of its differential 
 $f(0)(dz+\mu d\bar z)$ at 0 (the inequalities $|\mu|<1$, $f(0)\neq0$). 
 
 Recall that by (\ref{2.10}),    
$$f(x,t)=(Id-tU\circ\mu)^{-1}(1), \ \text{where}\ 
U=(\partial_{\bar z})^{-1}
 \partial_z.$$ The functions $f(x,t)$ are equal to 1, if 
 $\mu=0$. Let us show that they are $C^0$- close to 1 (and hence, 
$f(0,1)\neq0$), whenever $\mu$ is small enough with 
 derivatives up to order 3. Consider the operator functional 
 $\mathcal A(\mu)=(Id-tU\circ\mu)^{-1}$: its value being an operator 
 acting in $H^3(\td)$ (it is well-defined, see Proposition \ref{2.70}). 
 As it will be shown in the next paragraph, it  
 depends continuously on small functional parameter $\mu\in C^3(\td)$, 
 $max|\mu|<1$, in the 
 $H^3(\td)$ operator norm, and moreover, it has a bounded derivative in $\mu$. 
Therefore, if $||\mu||_{C^3}$ is small enough, then  each function $f(x,t)$ 
 is close to 1 in $H^3$ (thus, in $C^0$, by Sobolev embedding theorem). 

 Now for the proof of Proposition \ref{loc} it suffices to prove 
 the boundedness of the previous derivative $\mathcal A'(\mu)$. 
 For any $0<\delta'<1$ the $\mathcal A(\mu)$ is uniformly bounded in all 
 $\mu$ with $||\mu||_{C^3}<\delta'$ 
 (Proposition \ref{2.70}), so, we can apply the usual formula for 
 the derivative of the inverse operator: the derivative of $\mathcal A(\mu)$ 
 along a vector $h\in C^3(\td)$ is equal to 
 $$\nabla_h\mathcal A(\mu)=\mathcal A(\mu)\circ U\circ h\circ\mathcal A(\mu).$$
 To prove the boundedness of the derivative, we have to show that 
 the $H^3$- norm of the operator in the right-hand side of the previous formula 
 is no greater than some constant (depending on $\mu$) times $||h||_{C^3}$. 
 Indeed, the previous $H^3$ operator norm is no greater 
 than $||\mathcal A(\mu)||_{H^3}^2$ times the $H^3$- norm of the operator 
 of multiplication by the function $h$, the latter is no greater than 
 $||h||_{C^3}$ times some universal constant. This proves 
 the boundedness of the derivative. 
 Proposition \ref{loc} is proved. The proof of Theorem \ref{td} is completed.

 \section{Quasiconformal mappings. Proof of Theorem \ref{c}}
 
 \subsection{The plan of the proof of Theorem \ref{c}}
 
 We have already proved the statement of Theorem \ref{c} for a 
 $\ci$ double-periodic almost complex structure on $\cc$  
 (i.e., a lifting to the universal cover $\cc$ 
 of a $\ci$ complex structure on $\td$). In this case the diffeomorphism 
 $\cc\to\cc$ from the Theorem is the lifting to the universal covers 
 of the diffeomorphism of the tori given by Theorem \ref{td}.
 To prove Theorem \ref{c} in the general case (let $\sigma$ be a given 
 (may be measurable) bounded complex structure on $\cc$) we 
 consider a sequence $\sigma_n$ of $\ci$ double-periodic complex 
 structures on $\cc$ with growing periods and 
 uniformly bounded dilatations (say less than a fixed $K>0$)  
 that converge to $\sigma$ almost everywhere. 
 For each $\sigma_n$ there is a normalized quasiconformal 
 diffeomorphism $\Phi_n:\cc\to\cc$ 
 transforming $\sigma_n$ to the standard complex structure. 
 We show that the diffeomorphisms $\Phi_n$ converge 
 (uniformly on $\overline{\cc}$) to a homeomorphism (denoted $\Phi$). 
 We will prove that $\Phi$ is a quasiconformal homeomorphism sending 
 $\sigma$ to the standard complex structure 
 (see the end of the Subsection). The uniqueness of a latter homeomorphism 
 and its diffeomorphic property  
 on a smoothness domain of $\sigma$ will be proved in 3.4. 
 Its analytic dependence on 
 parameter (the Addendum to Theorem \ref{c}) will be proved in 3.5.  
 
 We prove the convergence of $\Phi_n$ by equicontinuity  
 of the normalized $K$- homeomorphisms: 
 
 \begin{lemma} \label{norm} [Ah2]. For any $K>0$ the normalized 
 $K$- homeomorphisms $\cc\to\cc$ (see Definition 1.1) 
 are equicontinuous with their inverses as mappings of the Riemann sphere. 
 \end{lemma}
 Lemma \ref{norm} (proved in  3.2) together with Arzela-Ascoli 
 theorem imply the following 
 \begin{corollary} \label{normc} For any $K>0$ each sequence of 
 normalized $K$- homeomorphisms $\cc\to\cc$ 
 contains a subsequence converging to a homeomorphism 
 $\overline{\cc}\to\overline{\cc}$ uniformly on $\overline{\cc}$. 
\end{corollary}

\begin{lemma} \label{limqc} [Ah2]. Let $K>0$, $U\subset\cc$ be a domain 
(that may be the whole $\cc$) $\Phi_n:U\to\Phi_n(U)\subset\cc$ be a sequence 
 of $K$- homeomorphisms converging uniformly on compact subsets to a 
 homeomorphism (denote $\Phi$ the limit). 
 Let $\sigma_n$ be the almost complex structures sent to the standard one by 
 $\Phi_n$. Let $\sigma_n$ converge almost everywhere (denote $\sigma$ 
 their limit). Then $\Phi$ is a $K$- homeomorphism  
 sending $\sigma$ to the standard complex structure. 
\end{lemma}

Lemma \ref{limqc} will be proved in Subsection 3.3 (using Lemma \ref{norm} and 
Corollary \ref{normc}).

\begin{proof} {\bf of existence in Theorem \ref{c} modulo Lemmas 
\ref{norm} and \ref{limqc}.} Let $\sigma_n$, $\sigma$, $K$, $\Phi_n$ be as at 
the beginning of the Section. Then $\Phi_n$ are $K$- 
diffeomorphisms. Passing to a subsequence, one can achieve that $\Phi_n$ 
converge to a homeomorphism (Corollary \ref{normc}, 
denote $\Phi$ the limit homeomorphism).  By Lemma \ref{limqc}, $\Phi$ is a $K$- 
homeomorphism 
transforming $\sigma$ to the standard complex structure. Theorem \ref{c} is proved. 
\end{proof}

\begin{remark} In the proof of the existence  in Theorem \ref{c} 
we had used only the statements of the previous 
Lemmas for $\ci$ diffeomorphisms. Their statements  
for general quasiconformal homeomorphisms will be used in the proof of 
the uniqueness in Theorem \ref{c} (Subsection 3.4). 
\end{remark}

\subsection{Normality. Proof of Lemma \ref{norm}}
The proof of Lemma \ref{norm} is based on the Gr\"otzsch inequality (the 
next Lemma) comparing moduli of $K$- homeomorphic complex annuli. To state it, 
let us firstly recall the following

\begin{definition} see [Ah2]. The {\it modulus} of an annulus 
$A=\{ r<|z|<1\}$ is $m(A)=-\frac1{2\pi}\ln r$. 
\end{definition}
\begin{remark} Consider the cylinder $\mathbb R\times S^1$ with the 
coordinates $(x,\phi)$, $S^1=\mathbb R\slash2\pi\mathbb Z$, and the standard 
complex structure, which is induced by the Euclidean metric 
$dx^2+d\phi^2$.  
\begin{equation}\text{For any}\ R>0 \ \text{put}\ A(R)=\{ 0<x<R\};\ \text{then}\ 
m(A(R))=\frac R{2\pi}.\label{ann}\end{equation}
The modulus of an annulus is invariant under conformal mappings [Ah2]. 
\end{remark} 
\begin{lemma} \label{grl} (Gr\"otzsch, see [Ah2]). Let $K>0$, 
$f:A_1\to A_2$ be a $K$- homeomorphism of complex annuli. Then 
\begin{equation} m(A_2)\geq K^{-1}m(A_1).\label{gri}\end{equation}
\end{lemma}
\begin{proof} For completeness of presentation, we give the classical proof 
of the Lemma. Firstly we prove the Lemma for a $K$- diffeomorphism; the general 
case is treated analogously (see the end of the 
proof). Let us consider that the annuli are drawn on the previous cylinder, 
say, $A_1=A(R_1)$, $A_2=A(R_2)$, then $m(A_i)=\frac{R_i}{2\pi}$, $i=1,2$, 
see (\ref{ann}).  Thus, it suffices to show that $R_2\geq K^{-1}R_1$. To do this, 
consider 
 the pullback (denoted $g$) to $A_1$ under $f$ of the Euclidean metric of $A_2$ 
 (denote $|\ |_g$ ($Area_g$) the corresponding norm of vector 
 fields on $A_1$ (respectively, the area), $Area$ being the Euclidean area). 
 One has $Area(A_i)=2\pi R_i$, $Area(A_2)=Area_g(A_1)$. We show that
 \begin{equation} Area_g(A_1)\geq K^{-1} Area(A_1).
 \label{area}\end{equation}
 This together with the previous formulas for the areas will prove the Lemma. 
 For the proof of (\ref{area}) we consider the family 
 $A(r)=\{ 0<x<r\}\subset A_1$ of subannuli in $A_1$, $r\leq R_1$, and prove 
 a lower bound of the derivative $(Area_g(A(r))'_r$. To do this, consider 
 the vector field $\frac{\partial}{\partial x}$ as the sum of its 
 component tangent to the circles $x=const$ and the $g$- orthogonal component 
 
 \begin{tabular}{ll}
\begin{minipage}{10em}
\begin{center}
\hskip-0.5cm
\epsfbox{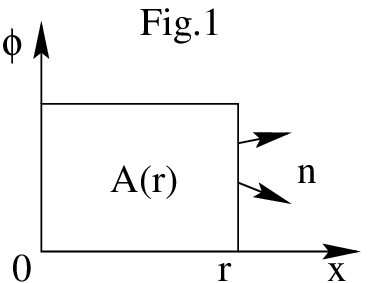}
\end{center}
\end{minipage}
\begin{minipage}{27em} 
(denote the latter component normal to the circles by $n$, see Fig.1). 
The vector field $n$ has the 
 same projection to the $x$- axis, as $\frac{\partial}{\partial x}$ and 
 its flow leaves invariant the fibration by circles $x=const$: its time $t$ 
 flow map transforms $A(r)$ to $A(r+t)$. Therefore, 
\end{minipage}
\end{tabular}
 \def\dpf{\frac{\partial}{\partial\phi}}
 \def\dpx{\frac{\partial}{\partial x}}
 \begin{equation}(Area_g(A(r)))'_r=\int_{x=r,\phi\in[0,2\pi]} 
 |\dpf|_g|n|_gd\phi.\label{arder}\end{equation}
  One has $|n|_g\geq K^{-1}|\dpf|_g$. Indeed, the $g$- norm $|\ |_g$ of a 
 vector tangent to $A_1$ 
 is equal to the standard Euclidean norm $|\ |$ of its image under $f$: 
 $|n|_g=|f_*n|$, $|\dpf|_g=|f_*\dpf|$. By definition,  
 $|\dpf|=1$, $|n|\geq|\frac{\partial}{\partial x}|=1=|\dpf|$. Therefore, 
 by the $K$- quasiconformality of $f$ (see, Definition 1.1), 
 $|n|_g=|f_*n|\geq K^{-1}|f_*\dpf|=K^{-1}|\dpf|_g$. Hence, 
the previous derivative is no less than 
$$K^{-1}\int_{x=r,\phi\in[0,2\pi]}|\dpf|_g^2d\phi\geq K^{-1}(2\pi)^{-1}
(\int_{\phi\in[0,2\pi]}|\dpf|_gd\phi)^2$$
(Cauchy-Bouniakovskii-Schwarz inequality). The latter integral is no less than 
$2\pi$. Indeed, it is equal to 
the length in the metric $g$ of the circle $x=r$, or in other terms, 
 the Euclidean length of its image under $f$, which 
is a closed curve in $A_2$ isotopic to a circle $x=const$. Therefore, 
$(Area_g(A(r)))'_r\geq2\pi K^{-1}$, thus,  
$Area_g(A_1)\geq 2\pi K^{-1}R_1=K^{-1}Area(A_1)$. 
This proves (\ref{area}) and the Lemma for a $K$- diffeomorphism $f$. 
In the case, when $f$ is a $K$- homeomorphism, thus just having  
local $L_2$ derivatives, the previous discussion remains valid: the 
previous integrals are well-defined for almost all $r$, since the subintegral 
expression $|\dpf|_g|n|_g$ in (\ref{arder}) is bounded from above by 
$||df(r,\phi)||^2$ times a constant depending on $K$. This follows 
from definition and the uniform boundedness of the Euclidean norm $|n|$:  
by definition,  $n$ is projected to the vector field 
$\frac{\partial}{\partial x}$ with unit norm; the  
angle between $n$ and a circle $x=const$ is bounded from below by a 
constant depending on $K$ (quasiconformality). Lemma \ref{grl} is proved. 
 \end{proof}
 
 \def\var{\varepsilon}

To prove Lemma \ref{norm}, we need to show that close points 
cannot be mapped to distant points under a normalized K- homeomorphism or its 
inverse. This is proved by comparing moduli of appropriate annuli with 
those of their images (using Lemma \ref{grl}). 

For the proof of Lemma \ref{norm} we recall the notion of the Poincar\'e metric 
[CG]. The Poincar\'e metric of the unit disc $|z|<1$ is 
$\frac{4|dz|^2}{(1-|z|^2)^2}$ (it is invariant under its conformal 
automorphisms). A Riemann surface is {\it hyperbolic}, 
if its universal covering is conformally equivalent to the unit disc 
(see Theorem \ref{pk}), e.g., any domain in $\cc$ whose complement contains 
more than one point. 
The Poincar\'e metric of a hyperbolic Riemann surface is the pushforward 
of the Poincar\'e metric of the unit disc under the universal covering. 

\begin{remark} \label{poim} (see [CG]). The Poincar\'e metric is well-defined, complete 
and decreasing: the Poincar\'e metric of 
a subdomain of a hyperbolic Riemann surface is greater than that of the ambient 
surface. The Poincar\'e metric of $\cc\setminus\{0,1\}$ is greater than 
its standard spherical metric times a constant. 
\end{remark}

In the proof of Lemma \ref{norm} we use the following relation of modulus 
of an annulus and its Poincar\'e metric, whose proof is a straightforward 
calculation.

\begin{proposition} \label{geod} see, e.g., [DH]. The modulus of an annulus is equal to 
 $\pi$ times the inverse of the length of its closed geodesic. 
\end{proposition}

Let us prove the equicontinuity of normalized $K$- homeomorphisms  
by contradiction. Suppose the contrary, i.e., 
there exist an $\var>0$, a sequence of normalized 
$K$- homeomorphisms $\Phi_n:\cc\to\cc$ and a sequence of pairs 
$x_n,y_n\in\cc$, $|x_n-y_n|\to0$, $|\Phi_n(x_n)-\Phi_n(y_n)|>\var$ 
(in the spherical metric of $\overline{\cc}$). Without loss of 
generality we assume that the sequence $x_n$ (and hence, $y_n$) converges 
(one can achieve this by passing to a subsequence, denote $x$ the limit). 
Then there is a sequence $A_n$ of annuli in 
$\overline{\mathbb C}\setminus\{0,1,\infty\}$ bounded by circles 
centered at $x$ and surrounding the pairs $x_n$, $y_n$: one of the circles 
is fixed, the other one contracts to $x$, as $n\to\infty$ see Fig.2a. 
 By definition, the annuli $A_n$ 
tend to once punctured disc, hence, 
$m(A_n)\to\infty$. The point $x$ may coincide with some of the 
three points 0, 1, $\infty$. Let us take two of the latters that 
are distinct from $x$ (say, let them be 0, 1). Then each annulus $A_n$ 
separates the pairs $(x_n,y_n)$ and $(0,1)$.  
By Lemma \ref{grl}, $m(\Phi_n(A_n))\to\infty$ as well. Hence, by  Proposition 
\ref{geod}, the 
lengths of the geodesics (denoted by $\gamma_n$) of the annuli $\Phi_n(A_n)$ 
in their Poincar\'e metrics tend to zero.
 But the latter lengths are greater than the lengths of $\gamma_n$ 
taken in the Poincar\'e metric of $\cc\setminus\{0,1\}$, and hence, 
also greater than their lengths  in the spherical metric times a constant independent 
from $n$ (by the previous Remark). Thus, each $\gamma_n$ separates the pairs 
$(\Phi_n(x_n), \Phi_n(y_n))$ and $(0,1)$ and is a closed curve with spherical 
length tending to 0. Hence, the spherical distance between 
$\Phi_n(x_n)$ and $\Phi_n(y_n)$ tends to 0 - a contradiction. 
\medskip

\begin{tabular}{ll}
\begin{minipage}{15em}
\begin{center}
\epsfbox{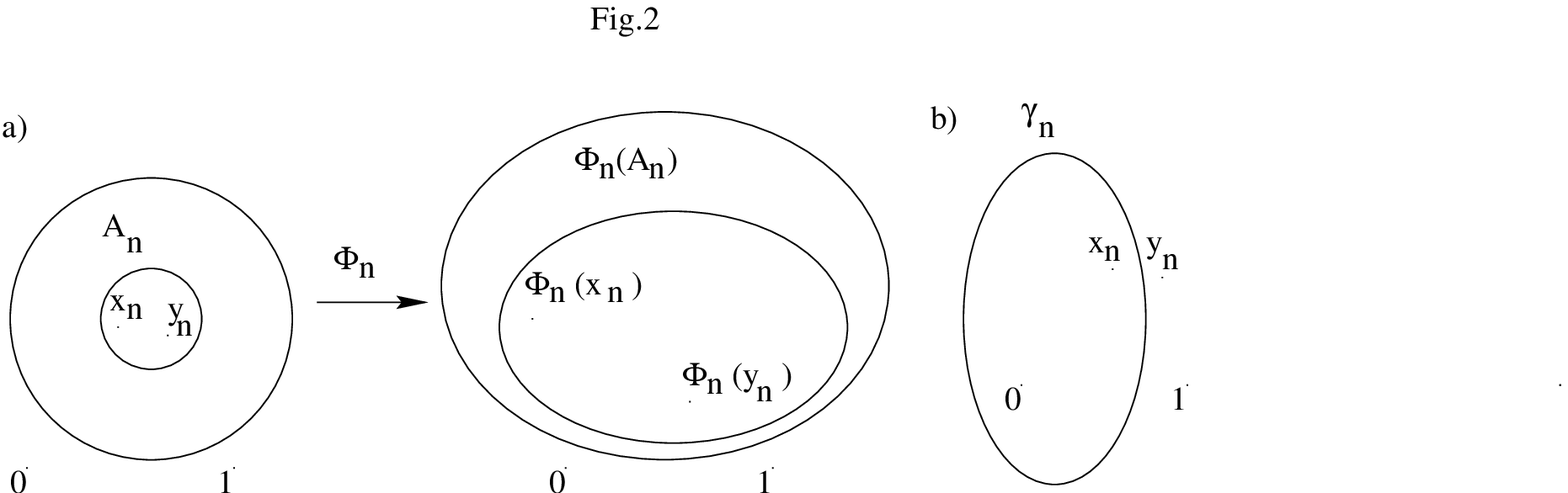}
\end{center}
\end{minipage} 
\end{tabular}
\medskip

Now let us prove that the inverses to the normalized $K$- homeomorphisms 
are also equicontinuous by contradiction, analogously to the previous 
discussion. Suppose the contrary:  
there exist an $\var>0$, a sequence of normalized 
$K$- homeomorphisms $\Phi_n:\cc\to\cc$ and a sequence of pairs 
$x_n,y_n\in\cc$, $|x_n-y_n|\to0$, $|\Phi_n^{-1}(x_n)-\Phi_n^{-1}(y_n)|>\var$ 
(in the spherical metric of $\overline{\cc}$). Without loss 
of generality we assume that all the 
sequences $x_n$, $y_n$, $\Phi_n^{-1}(x_n)$, $\Phi_n^{-1}(y_n)$ converge 
(one can achieve this by passing to a subsequence); denote their 
limits by $x$, $y$, $\tilde x$, $\tilde y$ respectively. 
By definition, $x=y$, $\tilde x\neq\tilde y$. 
Firstly consider the case, when $x\neq0,1,\infty$. Then $\tilde x,\tilde y
\neq0,1,\infty$ 
as well: otherwise $\Phi_n^{-1}(x_n)$, $\Phi_n^{-1}(y_n)$ would accumulate 
to $\{0,1,\infty\}$, while their $\Phi_n$- images $x_n$, $y_n$ would not - 
a contradiction to the equicontinuity of the $\Phi_n$'s (already proved). 
Fix an annulus $A$ separating the pair $0,\tilde x$ and the triple 
 $1,\tilde y,\infty$ (we assume that its closure is disjoint from 
$\Phi_n^{-1}(x_n)$, $\Phi_n^{-1}(y_n)$ for any $n$). 
Its images $\Phi_n(A)$ are disjoint from $0,1,x_n,y_n$ and 
have moduli bounded away from zero (Lemma \ref{grl}), 
and hence, closed geodesics (denoted $\gamma_n$, see Fig.2b) of 
uniformly bounded lengths. 
Thus, the lengths of $\gamma_n$  in the Poincar\'e metric of 
$\cc\setminus\{0,1,x_n,y_n\}$ are also uniformly bounded. On the other hand, 
$\gamma_n$ separates the pair $(0, x_n)$ and the triple $(1,y_n,\infty)$ 
for any $n$, see Fig.2b. The points 0, $x_n$ in the first pair are distant 
($x=\lim x_n\neq0$), thus, 
the spherical length of $\gamma_n$ is bounded from below. 
The points $x_n$, $y_n$, which are separated by $\gamma_n$, 
collide towards $x$, so, $\gamma_n$ comes arbitrarily close to $x_n$, as 
$n\to\infty$. This implies that $\gamma_n$ has length tending to infinity 
in the Poincar\'e metric of $\cc\setminus\{0,x_n\}$, and hence, in 
the Poincar\'e metric of $\cc\setminus\{0,1,x_n,y_n\}$. If   
$x_n$ does not move while $n$ changes, this  follows from the 
completeness of the Poincar\'e metric (Remark \ref{poim}). The case, when 
$x_n\not\equiv const$, is reduced to the previous one by 
applying the variable change $w=\frac z{x_n}$. This contradicts to the 
previous statement saying that the latter Poincar\'e length of $\gamma_n$ 
is uniformly bounded. 

Now let $x\in\{0,1,\infty\}$, say, $x=1$. 
Then $\tilde x,\tilde y\neq0,\infty$, as before, and at least one of 
$\tilde x\neq\tilde y$ (say, $\tilde x$) is distinct from 1. In these 
notations we repeat the previous argument. Lemma \ref{norm} is proved. 

\subsection{Quasiconformality and weak convergence. Proof of Lemma \ref{limqc}}
Let $\Phi_n$, $\sigma_n$, $\Phi$, $\sigma$ be as in Lemma \ref{limqc}. 
Recall that 
the dilatations of the $\sigma_n$' s are no greater than $K$, as are 
those of the $\Phi_n$' s, hence, the same is true for $\sigma$. Let us show 
that $\Phi$ is quasiconformal, more precisely: 1) has local $L_2$ derivatives 
that are weak $L_2$ limits of those of $\Phi_n$; 2) transforms $\sigma$ to the 
standard complex structure (and hence, is $K$- quasiconformal). 
This will prove Lemma \ref{limqc}.

For the proof of statement 1) we use the fact that the norms of 
the  differentials $d\Phi_n$ (in the spherical metric of $\overline{\cc}$) 
are uniformly bounded in each space $L_2(D)$, $D\Subset\cc$. 
Indeed, on each disc $D\Subset\cc$ 
$||d\Phi_n||^2_{L_2(D)}\leq K (Area(\Phi_n(D))$, which 
follows from definition and $K$- quasiconformality (the areas are taken  
in the spherical metric). The latter areas  
converge to $Area(\Phi(D))$, hence, they are uniformly bounded, and so are 
the previous $L_2$- norms. 

Thus, the derivatives are locally $L_2$- bounded, hence, passing to a 
subsequence one can achieve that they converge $L_2$- weakly. On the other 
hand, they converge to the 
derivative of $\Phi$ in sense of distributions. Therefore,  the latter is 
also $L_2$ locally and the convergence is $L_2$- weak. Statement 1) is proved. 

Let $\mu_n$, $\mu$ be the functions 
from (\ref{1.1}) defining the complex structures $\sigma_n$ and $\sigma$ 
respectively, thus, $d\Phi_n=f_n(dz+\mu_nd\bar z)$. By assumption, 
$|\mu_n|<1$, $\mu_n\to\mu$ almost everywhere. We have to show that 
$\frac{\partial \Phi}{\partial\bar z}=\mu\frac{\partial\Phi}{\partial z}$. 
 Indeed,  $f_n\to f=\frac{\partial\Phi}{\partial z}$, 
$f_n\mu_n\to\frac{\partial\Phi}{\partial\bar z}$  (both $L_2$ weakly), 
as $n\to\infty$. Since, $f_n$ are uniformly bounded in a local space $L_2$ 
and weakly converge, $\mu_n$ are uniformly bounded and converge 
almost everywhere, the weak limit of their product is the product 
$f\mu$ of their limits. This proves the previous partial differential 
equation on $\Phi$ together with statement 2) and Lemma \ref{limqc}. 
\subsection{Uniqueness, smoothness and group property}

Here we prove the uniqueness of the normalized 
quasiconformal homeomorphism from Theorem \ref{c} and the group 
and measure  
properties of quasiconformal mappings (Propositions \ref{group} 
and \ref{measure0}).  
The uniqueness follows from the local uniqueness up to composition with 
a conformal mapping and from normalizedness. The local uniqueness 
(together with the diffeomorphic property on a smoothness domain of 
the complex structure) are implied by the following 

\begin{proposition} \label{uniq} Let $D\subset\cc$ be a simply-connected 
domain, $\sigma$ be a bounded measurable  almost complex structure on $D$, 
$\Phi:D\to\Phi(D)\subset\cc$ be a quasiconformal homeomorphism transforming 
$\sigma$ to the standard complex structure. Then $\Phi$ is unique up to 
left composition with a conformal mapping. It is a $\ci$ diffeomorphism, 
if $\sigma$ is $\ci$. 
\end{proposition}

\begin{proof}  
Let $\mu:D\to\cc$ be the function defining the almost 
complex structure $\sigma$. 

{\bf Case $\mu\equiv0$.} Then $\frac{\partial \Phi}{\partial\bar z}=0$ and 
$\Phi$ has local $L_2$ derivatives. Let us show that $\Phi$ is conformal. 
Fix a $z_0\in D$ and put 
$U(z)=\int_{z_0}^z\Phi(\zeta)d\zeta$. We show that the function $U(z)$   
is well-defined (independent on the choice of path connecting $z_0$ 
to $z$). Then it is holomorphic by definition, hence, so is $\Phi(z)=
\frac{\partial U}{\partial z}$. It suffices to show that the integral 
of the form $\Phi dz$ along any Jordan curve is zero. Since the derivatives 
of $\Phi$ are locally $L_2$, we can apply the Stokes formula:  
the previous integral is equal to the integral of 
the differential $d(\Phi dz)$ over the domain bounded 
by the curve. But $d(\Phi dz)=\frac{\partial\Phi}{\partial\bar z}d\bar z dz=0$, 
so, it is zero.  

{\bf Case $\mu\in\ci$.} There exists at least one $\ci$ quasiconformal 
diffeomorphism $\Psi$ transforming $\sigma$ to the standard complex 
structure (Theorem \ref{td}, see also the discussion in 
Section 1.2). The composition $\Phi\circ\Psi^{-1}$ preserves 
the standard complex structure by definition and is quasiconformal: 
it has local $L_2$ derivatives, since so does $\Phi$ and $\Psi^{-1}$ is $\ci$.  
Therefore, it is conformal, as is proved above, hence, $\Phi$ is a $\ci$ 
diffeomorphism. 

{\bf Case $\mu$ is measurable.} 
Let $0<\delta<1$, $|\mu|<\delta$, $\mu_n$ be a sequence of $\ci$ functions, 
$|\mu_n|<\delta$, $\mu_n\to\mu$ almost everywhere (we extend 
$\mu_n$, $\mu$ to $\cc$ with the latter 
inequality and convergence). Consider the corresponding 
almost complex structures $\sigma_{\mu_n}$, see (\ref{1.1}), and 
the quasiconformal diffeomorphisms (denoted $\Phi_n$) from 
Theorem \ref{c} (the latters exist as is proved above). 
Passing to subsequence, one can assume that they converge uniformly 
on $\overline{\cc}$ (by Lemma \ref{norm}). Denote $\Psi$ their limit, 
which is a quasiconformal homeomorphism transforming the  extended 
complex structure $\sigma$ to the standard one (Lemma \ref{limqc}). 
It suffices to show that 
$\Phi\circ\Psi^{-1}:\Psi(D)\to\Phi(D)$ is a conformal homeomorphism. 
It is a homeomorphism, since so are $\Phi$ and $\Psi$, and preserves 
the standard complex structure, thus, if we show that it is 
quasiconformal, this will imply conformality 
(as is proved in the previous case $\mu\equiv0$). 
 To do this, 
consider the homeomorphisms $h_n=\Phi\circ\Phi_n^{-1}:\Phi_n(D)\to\Phi(D)$. 
They are quasiconformal homeomorphisms 
with uniformly bounded dilatations, as in the previous paragraph. They converge 
to $\Phi\circ\Psi^{-1}$ uniformly on compact subsets of $\Psi(D)$. The 
corresponding pullbacks of the standard complex structure converge 
to the latter almost everywhere, which follows from definition and 
convergence $\mu_n\to\mu$. Hence, by Lemma \ref{limqc}, the limit is 
quasiconformal. Proposition \ref{uniq} is proved. Theorem \ref{c} is proved.
\end{proof}

\begin{proof} {\bf of Proposition \ref{group}.} The statement of the 
Proposition is local: it suffices to show that compositions 
(inverses) of local $K$- homeomorphisms are $K^2$- (respectively, 
$K$-) quasiconformal. We prove this statement for composition 
(for inverse the proof is analogous): given domains $U,V,W\subset\cc$ 
and $K$- homeomorphisms $\Psi:U\to V$, $\Phi:V\to W$, let us show that 
$\Phi\circ\Psi$ is a $K^2$- homeomorphism. 
By Remark \ref{grr} the previous statement holds 
for diffeomorphisms, and in the general case the dilatation of the 
composition is no greater than $K^2$, thus, to prove the quasiconformality 
means to show that the composition has local $L_2$ derivatives. To do this, 
consider the pullback $\sigma(\Phi)$ of the standard complex 
structure under $\Phi$. Let us extend it to $\cc$ without increasing the 
dilatation and construct a sequence $\sigma(\Phi_n)$ of $\ci$ almost complex 
structures on $\cc$ converging to $\sigma(\Phi)$ 
almost everywhere with dilatations 
no greater than $K$. Let $\Phi_n:\cc\to\cc$ be the corresponding 
normalized quasiconformal homeomorphisms (which are $K$- homeomorphisms) from 
Theorem \ref{c}. They are $\ci$ diffeomorphisms as is proved above. 
By Lemma \ref{limqc}, they converge uniformly on $overline{\cc}$ to a 
$K$- homeomorphism $\wt{\Phi}:\cc\to\cc$ transforming 
 $\sigma(\Phi)$ to the standard complex structure. 
 By the previous Proposition, 
 $\wt{\Phi}=\Phi$ up to left composition with a conformal mapping. 
 Now the compositions $\Phi_n\circ\Psi$ are $K^2$- homeomorphisms 
 (since $\Phi_n$ is $\ci$) converging to $\wt{\Phi}\circ\Psi$, and the 
 corresponding pullbacks of the standard complex structure converge also. Hence, 
 by Lemma \ref{limqc}, the limit is quasiconformal. Since the limit coincides with 
 $\Phi\circ\Psi$ up to composition with a conformal mapping, the latter is 
 quasiconformal too. Proposition \ref{group} is proved.
\end{proof} 

\begin{proof} {\bf of Proposition \ref{measure0}.} The statement 
of Proposition \ref{measure0} is local and is reduced to the 
case of quasiconformal homeomorphisms $\overline{\cc}\to\overline{\cc}$, as 
Proposition \ref{group} proved above. We prove it 
by contradiction. Suppose the contrary: some quasiconformal 
homeomorphism $h$ of the Riemann sphere sends a zero measure set $S$ 
to a posivite measure set $h(S)$ (without loss of generality  we assume 
that $h$ fixes 0, 1 and $\infty$). 
Let $\sigma$ be the pull-back under $h$ 
of the standard complex structure  (it is 
well-defined almost everywhere). Then $h$ is the unique normalized 
quasiconformal homeomorphism transforming $\sigma$ to the standard 
complex structure. Let us change the standard structure in the image as 
follows: on $h(S)$ we change it 
to some constant nonstandard almost complex structure; on the rest we 
keep it standard. Denote  $\sigma'$ the almost 
complex structure thus obtained on the Riemann sphere in the image. 
By Theorem \ref{c}, there exists a unique normalized 
quasiconformal homeomorphism 
$H$ transforming $\sigma'$ to the standard complex  structure. 
One has $H\not\equiv Id$, since the set $h(S)$ has a positive measure. 
By definition and Proposition \ref{group}, the composition 
$H\circ h$, which is different from $h$, 
is a normalized quasiconformal homeomorphism transforming 
$\sigma$ to the standard complex structure. This contradicts the 
uniqueness of $h$. Proposition \ref{measure0} is proved. 
\end{proof}

\subsection{Analytic dependence on parameter. Proof of the Addendum}

{\bf Double-periodic case.} Consider a family of double-periodic $\ci$ 
almost complex 
structures $\sigma(t)$ on $\cc$ depending holomorphically 
on a complex parameter $t$ 
(this means that the corresponding function $\mu=\mu(z,t)$ from (\ref{1.1}) is 
holomorphic in $t$). We assume that the periods are fixed, thus, $\sigma(t)$ 
are the lifting to the universal cover $\cc$ of an analytic family of 
almost complex structures on the two-torus. 
Then the corresponding quasiconformal diffeomorphisms (denoted 
$\Phi_t$) from Theorem \ref{c} are holomorphic in $t$ as well. Indeed, 
their differentials are uniformizing differentials. Hence, for any $t$, 
$d\Phi_t=f_t(dz+\mu(z,t)d\bar z)$ up to multiplication by complex constant 
depending on $t$, where $f_t$ is given by formula (\ref{2.10}). The right-hand 
side of (\ref{2.10}) is analytic in the functional parameter $\mu$, hence, 
$f_t$ is holomorphic in $t$ and $z\mapsto\int_0^zf_t(dz+\mu(z,t)d\bar z)$ 
is a holomorphic family of diffeomorphisms of $\cc$. The family $\Phi_t$ is 
obtained from the latter by multiplication by a function in $t$ that 
makes the previous diffeomorphisms normalized (fixing 1), hence, the 
multiplier function (and thus, $\Phi_t$ as well) are also holomorphic in $t$. 
The Addendum is proved in the double-periodic case. 

{\bf General case.} 
Now consider arbitrary analytic family $\sigma(t)$ of bounded 
almost complex structures 
on $\cc$ depending on a complex parameter $t$ (we suppose that $t$ runs through 
the unit disc $D$). Let $\mu(z,t)$ be the corresponding functions, 
see (\ref{1.1}), which are holomorphic in $t$. Then there exists a 
$0<\delta<1$ such that $|\mu(z,0)|<\delta$ for any $z$. The  
corresponding mapping $M_z:t\mapsto\mu(z,t)$ is a holomorphic mapping 
$D\to D$ depending on $z$ in a measurable way such that $|M_z(0)|<\delta$. 
 (Recall that for a given $\delta<1$ the space of holomorphic mappings 
 $M:D\to D$ with $|M(0)|<\delta$ is compact, see [CG].) 
 Vice versa, for any $0<\delta<1$ each measurable collection of holomorphic 
 mappings $M_z:D\to D$ with $|M_z(0)|<\delta$ 
defines an analytic family of bounded almost complex structures; they 
are uniformly bounded when restricted to a smaller  
parameter disc $D_r=\{|t|<r\}$, $r<1$. Indeed, in the case, when $M_z(0)\equiv0$, 
$|M_z|_{D_r}<r$ (Schwarz Lemma); the general case is easily reduced to the 
previous one.  
 
Denote $\Phi_t$ the corresponding 
normalized quasiconformal homeomorphisms from Theorem \ref{c}.  
To prove the analyticity of $\Phi_t$ in $t$, we approximate $\sigma(t)$ 
(in the sense of convergence almost everywhere) 
by  analytic families $\sigma_n(t)$ of $\ci$ double-periodic almost 
complex structures depending holomorphically on the same parameter $t$ 
with growing periods $2n$, $2in$, 
$\sigma_n\to\sigma$. (For example, consider the restriction of $\sigma$ 
to the period square centered at 0 and take $\sigma_n$ to be its  
double-periodic  extension. Then approximate the new double-periodic  
family $M_z$ by a $\ci$ family of holomorphic mappings $D\to D$.) One can 
do this in such a way that $\sigma_n(t)|_{t\in D_r}$ be uniformly bounded.  
Denote $\Phi_{n,t}$ the normalized 
quasiconformal homeomorphisms transforming $\sigma_n(t)$ to the standard 
complex structure. They depend analytically on $t$, as is proved above. 
By Lemma \ref{limqc}, for any $t$, $z$, 
 $\Phi_{n,t}(z)\to\Phi_t(z)$, as $n\to\infty$. 
Thus, $\Phi_t(z)$ is a function in $t$ that 
is a limit of pointwise converging sequence of holomorphic functions. 
Let us prove that for any fixed $z$ the functions $\Phi_{n,t}(z)$ in $t\in D_r$ 
are bounded uniformly in $n$: then their limit $\Phi_t(z)$ is holomorphic. 
Indeed,  the almost complex structures $\sigma_n(t)|_{t\in D_r}$ 
are uniformly bounded. Therefore, the family $\Phi_{n,t}$ (depending 
on the two parameters $n$ and $t\in D_r$) together with their inverses 
is equicontinuous (Lemma \ref{norm}). Hence, the previous functions are 
uniformly bounded, so, their limit is holomorphic. The Addendum is proved. 

\section{Acknowledgements} 
In late 1990-ths \'E.Ghys stated a question concerning a foliated version of 
Theorem \ref{c} for linear foliations of tori. 
The proof of Theorem \ref{td} presented in the paper was obtained 
as a by-product of the author's solution to his question, and I wish to 
thank him.  
I wish also to thank him and  J.-P.Otal, \'E.Giroux for helpful discussions.
The research was supported by part 
by CRDF grant  RM1-2358-MO-02 and by 
Russian Foundation for Basic Research (RFFI) grant 02-02-00482. 

\section{References}

[Ab] Abikoff, W. Real analytic theory of Teichm\"uller space, - Lect. Notes in  
Math., 820, Springer-Verlag (1980). 

[Ah1] Ahlfors, L. Conformality with respect to Riemannian metrics. - 
Ann. Acad. Sci. Fenn. Ser. A. I. 1955 (1955), no. 206, 22 pp.
 
[Ah2] Ahlfors, L. Lectures on quasiconformal mappings, - Wadsworth (1987).  

[AhB] Ahlfors, L.; Bers, L. Riemann's mapping theorem for variable metrics, - 
Ann. of Math. (2) 72 (1960), 385-404. 

[Be] Bers, L. Riemann surfaces (mimeographed lecture notes), New York 
University (1957-58). 

[BeN] Bers, L.; Nirenberg, L.   
On a representation theorem for linear elliptic systems with discontinuous 
coefficients and its applications. 
Convegno Internazionale sulle Equazioni Lineari alle Derivate Parziali, Trieste, 1954, pp. 111--140. 
Edizioni Cremonese, Roma, 1955. 

[Bo] Boyarski\u\i, B. V. 
Generalized solutions of a system of differential equations of first order and of elliptic type with
discontinuous coefficients. (Russian) 
Mat. Sb. N.S. 43(85) 1957 451--503.

[CG] Carleson, L., Gamelin, Th.W. Complex Dynamics, - Springer-Verlag 1993. 

[Ch] Choquet-Bruhat, Y., \ de Witt-Morette, C., \ Dillard-Bleick, M. Analysis, 
Manifolds and Physics, - North-Holland, 1977.

[Chern] Chern, S.-S., 
An elementary proof of the existence of isothermal parameters on a surface. 
Proc. Amer. Math. Soc. 6 (1955), 771--782.

[DB] Douady, A.; Buff, X. Le th\'eor\`eme d'int\'egrabilit\'e 
des structures pr\`esque complexes. (French) [Integrability theorem
   for almost complex structures] The Mandelbrot set, theme and variations, 307--324, London Math. Soc. Lecture Note Ser., 274,
   Cambridge Univ. Press, Cambridge, 2000.

[DH] Douady, A.; Hubbard, J. 
A proof of Thurston's topological characterization of rational functions. 
Acta Math. 171 (1993), no. 2, 263--297. 

[Gl] Glutsyuk, A.; Simultaneous metric uniformization of foliations by 
Riemann surfaces. - To appear in Commentarii Mahtematici Helvetici. 

[Ko1] K\"obe, P. \"Uber die Uniformisierung beliebiger analytischer Kurven I. 
- Nachr. Acad. Wiss. G\"ottingen (1907), 177-190.

[Ko2] K\"obe, P. \"Uber die Uniformisierung beliebiger analytischer Kurven II. 
- Ibid. (1907), 633-669. 

[Korn] Korn, A., Zwei Anwendungen der Methode der sukzessiven Ann\"aherungen, - 
Schwarz Festschrift, Berlin (1919), pp. 215-229. 

[La] Lavrentiev, M.A., Sur une classe des repr\'esentations continues. - Mat. Sb., 
42 (1935), 407-434. 

[Licht] Lichtenstein, L., Zur Theorie der konformen Abbildungen; Konforme 
Abbildungen nicht-analytischer singularit\"atenfreier Fl\"achenst\"ucke auf 
ebene Gebiete, - Bull. Acad. Sci. Cracovie, (1916), 192-217. 

[M] Morrey, C. B., Jr. 
On the solutions of quasi-linear elliptic partial differential equations. - 
Trans. Amer. Math. Soc. 43 (1938), no. 1, 126--166.

[P] Poincar\'e, H. Sur l'uniformisation des fonctions analytiques, - 
Acta Math., 31 (1907), 1-64. 

[Vek]  Vekua, I. N. The problem of reduction to canonical form of differential forms of elliptic type and the generalized
   Cauchy-Riemann system. - (Russian) Dokl. Akad. Nauk SSSR (N.S.) 100, 
   (1955), 197--200.

\end{document}